\documentclass[reqno]{amsart}
\usepackage{amssymb,hyperref}
\usepackage[pagewise]{lineno}

\newtheorem{theorem}{Theorem}[section]

\newtheorem{definition}[theorem]{Definition}

\title{On existence and nonexistence of nonnegative solutions to seminlinear differential equation on Riemannian manifolds}
\author[Xu]{Fanheng Xu}
\address{School of Mathematical Sciences and LPMC, Nankai University, 300071 Tianjin, P. R. China}
\email{xufanheng@mail.nankai.edu.cn}


\begin{document}

\begin{abstract}
In this paper, we give a clear cut relation between the the volume growth $V(r)$ and the existence of nonnegative solutions to parabolic semilinear problem
\begin{align}\tag{*}\label{*}
\left\{
\begin{array}{ll}
\Delta u - \partial_t u + u^p = 0, \\
{u(x,0)= {u_0(x)}},
\end{array}
\right.
\end{align}
on a large class of Riemannian manifolds.
We prove that for parameter $p>1$, if 
\begin{align*}
\int^{+\infty} \frac{t}{V(t)^{p-1}} dt = \infty
\end{align*}
then (\ref{*}) has no nonnegative solution.
If 
\begin{align*}
\int^{+\infty} \frac{t}{V(t)^{p-1}} dt < \infty
\end{align*}
then (\ref{*}) has positive solutions for small $u_0$.
\end{abstract}
\maketitle

\section{introduction}
Let $(M,g)$ be a noncompact complete connected $n$-dimensional Riemannian manifold and consider nonnegative solutions to differential equations
\begin{align}\label{para_def_eq}
\left\{
\begin{array}{ll}
\Delta u - \partial_t u + u^p = 0 \quad \mbox{in} \ M \times (0, \infty ), \\
{u(x,0)= {u_0(x)}} \quad \mbox{in}\ M,
\end{array}
\right.
\end{align}
where parameter $p>1$, $\Delta$ is the Laplace-Bettrami operator of $M$, and $u_0(x) \geqslant 0$ not identically zero.

Fujita in \cite{Fujita} proved the following results in the case of $M=\mathbb{R}^n$:
\begin{itemize}
\item[(a)] if $1<p<1+\frac{2}{n} $ and $u_0>0$, equation (\ref{para_def_eq}) has no solution;
\item[(b)] if $p>1+\frac{2}{n} $ and $u_0>0$ is smaller than a small Gaussian, then equation (\ref{para_def_eq}) possesses positive solutions.
\end{itemize} 
Later, several authors 
\cite{Aronson_Weinberger} 
\cite{Hayakawa} 
\cite{Kobayashi_Sirao} 
showed that the critical case for problem (\ref{para_def_eq}), that is $p=1+\frac{2}{n} $, belongs to blow-up case (a). 

In \cite{Zhang_duke}, Zhang proved several nonexistence results for semilinear and quasilinear parabolic problems on manifolds provided the following setting:
\begin{itemize}
\item[(i)] $V_{x_0}(r) \leqslant C r^\alpha$;
\item[(ii)] Condition (G) (see Definition \ref{assumption_G}).
\end{itemize}
Especially, he proved that: Assume  conditions (i) and (ii) on manifold are satisfied, and $\alpha\geq1$. If $1<p\leq 1+\frac{2}{\alpha}$, then problem (\ref{para_def_eq}) possesses no global positive solution.

The celebrated idea of studying the nonnegative solutions in terms of the volume of the geodesic ball was due to Cheng and Yau \cite{Cheng_Yau}. 
They proved that if the inequality
\begin{align*}
V_{x_0}(r)\leqslant Cr^2
\end{align*}
holds for some point $x_0\in M$ with all large enough r, then any positive solution to $\Delta u\leqslant 0$ is identically constant, here $V_{x_0}(r)$ is the volume of geodesic ball of radius $r$ centered at $x_0 \in M$.

Recently, this idea was used and developed for solutions to both elliptic equations (see, e.g. 
\cite{Grigoryan_Kondratiev, Grigoryan_Sun, Sun_jmaa, Wang_Xiao})
and parabolic equations (see 
\cite{Mastrolia_Monticelli_Punzo_parabolic}).
Particularly in \cite{Mastrolia_Monticelli_Punzo_parabolic}, 
Mastrolia, Monticelli and Punzo proved that if
\begin{align}\label{vc_mmp}
V_{x_0}(r) \leqslant C r^\frac{2}{p-1} (\ln r)^\frac{1}{p-1},
\end{align}
then (\ref{para_def_eq}) has no weak supersolutions. 

\begin{definition}
We say the volume condition
\begin{align}\label{f2bc7fb1}
\int^{+\infty} \frac{t}{V(t)^{p-1}} dt =(or<) \infty
\end{align}
 is satisfied if there exists constant $r_0>0$ and $x_0 \in M$ such that
\begin{align*}
\int_{r_0}^{+\infty} \frac{t}{V_{x_0}(t)^{p-1}} dt = (or<) \infty
\end{align*}
\end{definition}
Noting that the choice of $r_0$ and $x_0$ doesn't affect this condition, we can omit them in (\ref{f2bc7fb1}) and specify proper value if needed.

It's easily see that, if $V_{x_0}(r)$ satisfy (\ref{vc_mmp}) or the following inequality 
\begin{align*}
V_{x_0}(r) \leqslant C r^{\alpha_1} (\ln r)^{\alpha_2} (\ln \ln r)^{\alpha_3} \times \cdots \times
\mathop{(\underbrace {\ln \ln \cdots \ln r})^{\alpha_k}}
\limits_{\text{k-1\quad}}, \quad k\geqslant 3
\end{align*}
with parameter
\begin{align*}
\alpha_1 = \frac{2}{p-1},\quad \alpha_i = \frac{1}{p-1},\ (i=2,3,\cdots,k),
\end{align*}
then $V_{x_0}(r)$ satisfy
\begin{align}\label{vc_nonexistence}
 \int^{+\infty} \frac{t}{V(t)^{p-1}} dt = \infty.
\end{align} 
In this paper, we show that volume condition (\ref{vc_nonexistence}) is indeed the sufficient condition to the nonexistence of solutions to equation (\ref{para_def_eq}) on a large class of manifolds. 
Moreover, if we assume the initial value $u_0$ is small, then (\ref{vc_nonexistence}) is also necessary.

Let us introduce our setting.
Unless specified, let $M$ be a connected geodesically complete noncompact Riemannian manifold and satisfy condition (G) and (H).

\begin{definition}[Condition (G)]\label{assumption_G}
We say condition (G) is satisfied if 
there exists a constant $C_0$ such that
\begin{align*}
\frac{\partial \log g^{1/2}}{\partial r} \leqslant \frac{C_0}{r}, 
\end{align*}
here $g^{1/2}$ is the volume density of the manifold;
\end{definition}

If the Ricci curvature of $M$ is nonnegative, it is well known that condition (G) hold, in fact, $(\partial \log g^{1/2})/\partial r \leqslant 0$.

\begin{definition}[Condition (H)]
We say condition (H) is satisfied if $P_t(x,y)$,
the smallest fundamental solution of the classical heat equation
\begin{align*}
\frac{\partial}{\partial t} u -\Delta u =0, \text{ in } M
\end{align*}
exists and satisfies the following estimate 
\begin{align}\label{hk_upper}
P_t(x,y)\leqslant \frac{C_1}{V_x(\sqrt{t})},
\end{align}
for all $t>0$ and almost all $x, y \in M$.
\end{definition}
We know that $P_t(x,y)$ is called heat kernel of $\Delta$ and has the following properties
\begin{itemize}
\item Symmetry: $P_t(x,y)= P_t(y,x)$, for all $x,y\in M, t>0$.
\item $P_t(x,y) \geqslant 0$, for all $x$, $y \in M$ and $t>0$, and
\begin{align}\label{hk_markov}
\int_{M} P_t(x,y)dy \leqslant 1,\quad\mbox{for all $x\in M$ and $t>0$}.
\end{align}
\item The semigroup identity: for all $x$, $y \in M$ and $t$, $s>0$,
\begin{align}\label{hk_semigroup}
P_{t+s}(x,y) = \int_{M} P_t(x,z)P_s(z,y)dz.
\end{align}
\end{itemize} 

If the manifold has nonnegative Ricci curvature, 
by a famous result of Li and Yau in \cite{Li_Yau}, we have
\begin{align*}
P_t(x,y)\asymp \frac{C}{V_x(\sqrt{t})}\exp\left(-\frac{d^2(x,y)}{ct}\right),
\end{align*}
where the sign $\asymp$ means that both $\leqslant$ and $\geqslant$ are true but with different values of
$C, c$. It follows that (H) is satisfied if the manifold has nonnegative Ricci curvature.

The main results of the paper are the following.
\begin{theorem}\label{thm_para}
(a). If 
\begin{align*}
\int^{+\infty} \frac{t}{V(t)^{p-1}} dt = \infty
\end{align*}
then equation (\ref{para_def_eq}) has no nonnegative solutions.

(b). If 
\begin{align*}
\int^{+\infty} \frac{t}{V(t)^{p-1}} dt < \infty
\end{align*}

and $u_0$ is smaller than a small Gaussian, then (\ref{para_def_eq}) has positive solutions. 
\end{theorem}

In the part (a) of Theorem \ref{thm_para}, the ``solution'' can also be understood in a weak sense. 
Our proof are based on a method originating from 
\cite{Zhang_duke} (also see \cite{Zhang_jmaa}) and the proof of Theorem 11.14 in \cite{Grigoryan_book},
that uses the Laplacian of the distance function and a carefully chosen test function. The arguments of part (b) are purely functional analytic. 

We denote by $W_{c}^2$ the subspace of $W_{loc}^2$ of functions with compact support,
and by $B_x(r)$ the geodesic ball of radius $r$ centered at $x \in M$.
$V_x(r)$ is the volume of $B_x(r)$.
In all the proof, since the geodesic ball always centered at a fixed point $x_0 \in M$,
we may use the abbreviation $B\equiv B_{x_0}$ and $V\equiv V_{x_0}$.
The symbol $C$, $C_0$, $C_1$, $\cdots$ denote positive constants whose values are unimportant and may vary at different occurrences.

\section{Proof of Theorem \ref{thm_para} }
\subsection*{Proof of part (a)}
Suppose $u$ is a solution of (\ref{para_def_eq}) and $\psi\geqslant 0$ is a function in $W_{c}^{2}(M\times [0, \infty ))$.
Multiplying $\psi$ to both side of (\ref{para_def_eq}) and integrating by part, we have
\begin{align}\label{para_def_weak_weak}
\int_0^\infty \int_{M} u^p\psi dx dt &= -\int_0^\infty \int_{M} u\Delta \psi dx dt - \int_0^\infty \int_{M} u \partial_t\psi dx dt - \int_Mu_0\psi dx \notag\\
&\leqslant -\int_0^\infty \int_{M} u\Delta \psi dx dt - \int_0^\infty \int_{M} u \partial_t\psi dx dt
\end{align}
We use the following test function:
\begin{align*}
\psi(x,t) = \varphi^q(r(x),t)
\end{align*}
where $q=p/(p-1)$, and $r(x) = d(x_0,x)$ for some fixed point $x_0 \in M$.
It follows that
\begin{align*}
\Delta \psi = q\varphi^{q-1}\Delta\varphi+q(q-1)\varphi^{q-2}\left| {\nabla \varphi } \right|^2,
\end{align*}
\begin{align*}
\partial_t\psi = q\psi^{q-1}\partial_t\psi
\end{align*}
Subtituding above into (\ref{para_def_weak_weak}) we have
\begin{align}\label{para_en_taro_adun}
\int_0^\infty \int_{M} u^p \varphi^q dx dt \leqslant 
C\left(\int_{M} u \varphi^{q-1}(-\Delta\varphi) dx + 
\int_{M} u \varphi^{q-1}(-\partial_t\varphi) dx \right)
\end{align}
Let $h\in W_c^2[0, \infty )$ be a function satisying
\begin{center}
$h(r)=1$, $r \in \left[ 0, 1 \right]$; $h(r)=0$, $r \in \left[ 2, \infty \right)$;
$-C_1 \leqslant h' \leqslant 0$, $\left| h'' \right| \leqslant C_1$, $r \in (1,2)$.
\end{center}
Fix a finite increasing sequence $\{r_k\}$, $({k=0,1,\cdots,i})$.
Let $Q_k = B(r_k)\times [0,r_k^2)$. 
Define $\varphi$ by
\begin{align*}
\varphi (x,t) = \left\{ {\begin{array}{*{20}{l}}
 1, & (x,t) \in Q_0,\\ 
 a \frac{(r_k-r_{k-1})^2}{V(r_k)^{p-1}}h(\frac{r(x)}{r_{k-1}})h(\frac{t}{r_{k-1}^2}) + T_k, &
(x,t) \in Q_k \backslash Q_{k-1}, k=1, \cdots, i\\ 
0, & (x,t) \in M\backslash Q_i
\end{array}} \right.
\end{align*}
where
\begin{align*}
a=\left( \sum\limits_{k=1}^i \frac{(r_k-r_{k-1})^{2}}{V(r_k)^{p-1}}\right)^{-1}
\end{align*}
and
\begin{align*}
T_k = \left\{ 
\begin{array}{*{20}{l}}
a\sum\limits_{j=k+1}^i \frac{(r_j-r_{j-1})^{2p}}{V(r_j)^{p-1}}, & 
k=1, \cdots, i-1\\
0, & k=i
\end{array}
\right.
\end{align*}
Clearly $\varphi \in W_c^{2}(M\times[0,\infty))$. 
Noting $\varphi(\cdot,t)$ is radial on $M$, we have when $(x,t) \in Q_k\backslash Q_{k-1}$, 
\begin{align}\label{para_Delta_varphi}
\Delta \varphi = \frac{\partial^2 \varphi(r,t)}{\partial r^2}
 +\left(\frac{n-1}{r}+\frac{\partial \log g^{1/2}}{\partial r} \right)
\frac{\partial \varphi(r,t)}{\partial r}
\end{align}
Specifying $\{r_i\}$ to be a geometric sequence with $r_k=2r_{k-1}$, we have 
\begin{align}\label{para_varphi'r}
-Ca\frac{(r_k-r_{k-1})}{V(r_k)^{p-1}} \leqslant 
\frac{\partial \varphi(r,t)}{\partial r} \leqslant 0,
\end{align}
and 
\begin{align}\label{para_varphi''r}
\left| \frac{\partial^2 \varphi(r,t)}{\partial r^2}\right| \leqslant \frac{Ca}{V(r_k)^{p-1}}.
\end{align}

Combining (\ref{para_Delta_varphi})， (\ref{para_varphi'r}), (\ref{para_varphi''r}) and condition (G), we have
\begin{align}\label{para_delta_varphi}
- \Delta\varphi \leqslant \frac{Ca}{V(r_k)^{p-1}},
\quad (x,t) \in Q_k\backslash Q_{k-1}.
\end{align}
Also by the definition of $\varphi$, we have 
\begin{align}\label{para_partial_t_varphi}
- \partial_t \varphi \leqslant \frac{Ca}{V(r_k)^{p-1}},
\quad (x,t) \in Q_k\backslash Q_{k-1}.
\end{align}
Subtituding (\ref{para_delta_varphi}) and (\ref{para_partial_t_varphi}) into (\ref{para_en_taro_adun}), we have
\begin{align*}
\int_0^\infty \int_{M} u^p \varphi^q dx dt \leqslant 
Ca \int_0^\infty \int_{M} \sum\limits_{k=1}^i \frac{u \varphi^{q-1}}{V(r_k)^{p-1}} \chi_{\{ Q_k\backslash Q_{k-1} \}} dx dt 
\end{align*}
Using the H\"{o}lder's inequality, we obtain 
\begin{align*}
& \int_0^\infty \int_{M} u^p \varphi^q dx dt \notag\\
&\leqslant 
Ca \left( \int_0^\infty \int_{M} u^p \varphi^q dx dt \right)^{1/p}
\left( \int_0^\infty \int_{M} \left( \sum\limits_{k=1}^i \frac{\chi_{\{ Q_k\backslash Q_{k-1} \}}}{V(r_k)^{p-1}} \right)^q dx dt \right)^{1/q} \notag \\
&=
Ca \left( \int_0^\infty \int_{M} u^p \varphi^q dx dt \right)^{1/p}
\left( \sum\limits_{k=1}^i \iint_{Q_k\backslash Q_{k-1}} \frac{1}{V(r_k)^{p}} dx dt \right)^{1/q}
\end{align*}
It follows that
\begin{align}\label{para_I_le_a}
\left( \int_0^\infty \int_{M} u^p \varphi^q dx dt \right)^{(p-1)/p}
&\leqslant Ca \left( \sum\limits_{k=1}^i \iint_{Q_k\backslash Q_{k-1}} \frac{1}{V(r_k)^{p}} dx dt \right)^{1/q} \notag\\
&\leqslant Ca \left( \sum\limits_{k=1}^i 
\frac{r_k^2}{V(r_k)^{p-1}} \right)^{1/q} \notag\\
&\leqslant Ca \left( \sum\limits_{k=1}^i \frac{(r_k-r_{k-1})^{2}}{V(r_k)^{p-1}}\right)^{1/q} \notag \\
&= Ca^{(q-1)/q}
\end{align}
On the other hands
\begin{align}\label{para_a^-1}
a^{-1}&=\sum\limits_{k=1}^{i} \frac{(r_k-r_{k-1})^{2}}{V(r_k)^{p-1}} \notag\\
&= C\sum\limits_{k=1}^{i} \frac{r_k^{2}-r_{k-1}^{2}}{V(r_k)^{p-1}} \notag\\
&\geqslant
C\sum\limits_{k=2}^{i} \int_{r_{k-1}}^{r_{k}} \frac{r}{V(r)^{p-1}}dx \notag\\
&= 
C\int_{2r_0}^{r_i} \frac{r}{V(r)^{p-1}}dx
\end{align}
Using the volume condition
\begin{align*}
\int^\infty \frac{r}{V(r)^{p-1}} dx = \infty,
\end{align*}
for every $r_0 > 0$, we have form (\ref{para_a^-1}), if $i \rightarrow +\infty$, then $a \rightarrow 0$. Thus we can choose sufficient large $i$, such that
\begin{align*}
 a \leqslant \frac{1}{r_0} 
\end{align*}
Subtituding above into (\ref{para_I_le_a}) , we have
\begin{align*}
\int_0^{r_0} \int_{B(r_0)} u^p dx dt
\leqslant
\frac{C}{r_0^{q-1}}
\end{align*}
Let $r_0 \rightarrow +\infty$, we obtain
\begin{align*}
\int_0^\infty \int_{M} u^p dx dt \rightarrow 0
\end{align*}
which implies that 
\begin{align*}
u \equiv 0.
\end{align*}
But $u_0\geq0$ is not identically zero, by Maximum principle, we know that $u$ is positive almost everywhere, which leads a contradiction, and then there's no nonnegative solution to problem (\ref{para_def_eq}).

\subsection*{Proof of part (b)}
Fix a point $x_0\in M$.
Let us define the operator
\begin{align}\label{T_u}
Tu(x,t)=\int_{M} P_t(x,y)u_0(y)dy + \int_0^t \int_{M} P_{t-s}(x,y)u^p(y,s) dx(y)ds.
\end{align}
acting on the space $S_M$ defined by
\begin{align}\label{S_M}
S_M=\left\{u\in L^{\infty}(M\times[0,\infty))|0 \leqslant u(x,t) \leqslant \lambda P_{t+\delta}(x,x_0)\right\}.
\end{align}
where $\lambda >0$ is a constant to be chosen later, and $\delta>1$ is a fixed constant.
It is easy to see that $S_M$ is a close set of $L^{\infty}(M\times[0,\infty))$.

Let $u_0$ satisfy
\begin{align}\label{assumption: u_0}
0\leqslant u_0(x)\leqslant \frac{\lambda}{2}P_{\delta}(x,x_0).
\end{align}
Now let us show $TS_M\subset S_{M}$.
Using (\ref{hk_semigroup}) and (\ref{assumption: u_0}), we have
\begin{align}\label{eq: Tu_part_1}
\int_{M}P_t(x,y)u_0(y)dx(y)&\leqslant \frac{\lambda}{2}\int_{M}P_t(x,y)P_{\delta}(y,x_0)dx(y)\notag\\
&=\frac{\lambda}{2}P_{t+\delta}(x,x_0).
\end{align}
Using (\ref{hk_upper}), (\ref{hk_semigroup}) and (\ref{S_M}), we have
\begin{align}\label{eq: en_taro_adun}
&\int_0^t \int_{M} P_{t-s}(x,y)u^p(y,s) dx(y)ds \notag\\
&\leqslant \lambda^p\int_0^t \int_{M} P_{t-s}(x,y)P_{s+\delta}^p(y,x_0) dx(y)ds \notag\\
&\leqslant \lambda^pC_3^{p-1}\int_0^t \frac{1}{V(\sqrt{s+\delta})^{p-1}} \int_{M} P_{t-s}(x,y)P_{s+\delta}(y,x_0) dx(y)ds \notag\\
&\leqslant \lambda^p C_3^{p-1} P_{t+\delta}(x,x_0) \int_0^t \frac{1}{V(\sqrt{s+\delta})^{p-1}} ds
\end{align}
By our volume condition
\begin{align*}
\int^{\infty} \frac{s}{V(s)^{p-1}} ds < \infty,
\end{align*}
there exists a constant $C_4$ such that 
\begin{align}\label{eq: int_VB_leq_C3}
\int_0^t \frac{1}{V(\sqrt{s+\delta})^{p-1}} ds
&= \int_\delta^{t+\delta} \frac{1}{V(\sqrt{s})^{p-1}} ds \notag \\
&= C\int_{\sqrt{\delta}}^{\sqrt{t+\delta}} \frac{s}{V(s)^{p-1}} ds \notag \\
&\leqslant C_4
\end{align}
Subtituding above into (\ref{eq: en_taro_adun}), we obtain, for small enough $\lambda$,
\begin{align}\label{eq: Tu_part_2}
\int_0^t \int_{M} P_{t-s}(x,y)u^p(y,s) dx(y)ds &\leqslant \lambda^p C_1^{p-1}C_4  P_{t+\delta}(x,x_0) \notag\\
&\leqslant \frac{\lambda}{2} P_{t+\delta}(x,x_0).
\end{align}
Combining (\ref{T_u}),(\ref{eq: Tu_part_1}) with (\ref{eq: Tu_part_2}), we obtain
\begin{align*}
0 \leqslant Tu \leqslant \lambda P_{t+\delta}(x,x_0).
\end{align*}
Hence
\begin{align*}
TS_{M} \subset S_M.
\end{align*}

Now let us show that $T$ is a contraction map. For $u_1$, $u_2\in S_M$, we have
\begin{align}
\left|Tu_1(x,t)-Tu_2(x,t)\right| \leq
\int_0^t \int_{M} P_{t-s}(x,y)\left|u_1^p(y,s)-u_2^p(y,s)\right| dx(y)ds
\end{align}
Noting that
\begin{align*}
\left|u_1^p(y,s)-u_2^p(y,s)\right| \leqslant p \max\{u_1^{p-1}(y,s),u_2^{p-1}(y,s)\}\left|u_1(y,s)-u_2(y,s)\right|
\end{align*}
and combining with (\ref{hk_upper}), (\ref{hk_markov}), (\ref{S_M}) and (\ref{eq: int_VB_leq_C3}), we obtain
that
\begin{align*}
& \left|Tu_1(x,t)-Tu_2(x,t)\right| \notag\\
&\leqslant p\lambda^{p-1}\left\| u_1-u_2 \right\|_{L^\infty}\int_0^t \int_{M} P_{t-s}(x,y)P_{s+\delta}^{p-1}(y,x_0) dyds \notag\\
&\leqslant p\lambda^{p-1}C_1^{p-1}\left\| u_1-u_2 \right\|_{L^\infty}
\int_0^t \frac{1}{V(\sqrt{s+\delta})^{p-1}}
\int_{M} P_{t-s}(x,y) dyds \notag\\
&= p\lambda^{p-1}C_1^{p-1}\left\| u_1-u_2 \right\|_{L^\infty}
\int_0^t \frac{1}{V(\sqrt{s+\delta})^{p-1}} ds \notag\\
&\leqslant p\lambda^{p-1}C_1^{p-1}C_4\left\| u_1-u_2 \right\|_{L^\infty}
\end{align*}

Readjusting $\lambda$ small enough such that $p\lambda^{p-1}C_1^{p-1}C_4<1$, we obtain that $T$ is a contraction map.
Applying fixed point theorem, we derive that
there exists a fixed point $u\in S_M$ satisfying
\begin{align}\label{fix}
u(x,t)=\int_{M} P_t(x,y)u_0(y)dy + \int_0^t \int_{M} P_{t-s}(x,y)u^p(y,s) dyds.
\end{align}
Since $u_0\gneqq0$, then $u$ is positive on $M$.
Since $u_0, u\in L^2(M)$, 
by a standard argument of regularity (see Theorem 7.6 and 7.7, \cite{Grigoryan_book}), 
we know the integrals in (\ref{fix}) are both smooth on $M\times(0,\infty)$, hence
we obtain that $u$ is a global positive solution of problem (\ref{para_def_eq}).

\end{document}